\documentclass[12pt]{article}    % Specifies the document style.
\usepackage{amsmath}
\usepackage{amssymb}
\usepackage{amsmath,amssymb,amsfonts,amsthm,latexsym}
\usepackage{booktabs}
\usepackage{array} 
\usepackage{lscape}
\usepackage{graphicx}
\usepackage[all]{xy}
\usepackage{epsfig}
\usepackage{colordvi}
\usepackage{latexsym}
\usepackage{authblk}
\def\mysection{\setcounter{equation}{0}\section}
\newtheorem{Lemma}{Lemma}[section]
\newtheorem{Theorem}[Lemma]{Theorem}
\newtheorem{Proposition}[Lemma]{Proposition}

\theoremstyle{definition}
\newtheorem{defn}{Definition}[section]

\newcommand{\beq}{\begin{equation}}
\newcommand{\eeq}{\end{equation}}
\newcommand{\beqr}{\begin{eqnarray}}
\newcommand{\eeqr}{\end{eqnarray}}

\begin{document}
	\title{\Large {{Topogenous structures and related families of morphisms} }}
	\author[$*,\;a$]{Minani Iragi}
	\author[$b$]{David Holgate}
	
	\affil[$a$]{{\footnotesize Institute of Mathematics, Faculty of Mechanical Engineering, Brno
			University of Technology\\ 
			Technická 2, 616 69 Brno, Czech Republic.}}
	
	\affil[$b$]{{\footnotesize Department of Mathematics and Applied Mathematics\\ University of the Western Cape,
			Bellville 7535, South Africa.
   %and Institute of Mathematics, 
   %Faculty of Mechanical Engineering, Brno
		%	University of Technology\\ 
			%Technická 2, 616 69 Brno, Czech Republic.
   }}
	
	\date{}
	\maketitle
	\def\thefootnote{\fnsymbol{footnote}}
	\setcounter{footnote}{0}
	
	\footnotetext{ 
		E-mail addresses:  $^{a}$Minani.Iragi@vutbr.cz   $^{b}$dholgate@uwc.ac.za\\
		$^{*}$Corresponding author.\\
	The first author acknowledges the support from the Brno University of Technology (BUT) under the project MeMoV II no. CZ.02.2.69/0.0/0.0/18-053/0016962. 
		%The Second author acknowledges the National Research Foundation of South Africa and the support from the European Social Fund and the state budget of the Czech Republic in the project no. CZ.02.2.69/0.0/0.0/16-027/0008371.
  }
	\begin{abstract}
	In a category $\mathcal{C}$ with a proper $(\mathcal{E}, \mathcal{M})$-factorization system, we study the notions of strict, co-strict, initial and final morphisms with respect to a topogenous order. Besides showing that they allow simultaneous study of four classes of morphisms obtained separately with respect to closure, interior and neighbourhood operators, the initial and final morphisms lead us to the study of topogenous structures induced by pointed and co-pointed endofunctors. We also lift the topogenous structures along an $\mathcal{M}$-fibration. This permits one to obtain the lifting of interior and neighbourhood operators along an $\mathcal{M}$-fibration and includes the lifting of closure operators found in the literature. A number of examples presented at the end of the paper demonstrates our results.
\end{abstract}
\\

\mysection{Introduction}
{\indent  The notion of categorical closure operators (\cite{dikranjan1987closure}) introduced by Dikranjan and Giuli was the point of departure for studying topological structures in arbitrary categories. This approach eventually motivated the introduction of categorical interior (\cite{MR1810290}) and neighbourhood (\cite{MR2838385}) operators. While the categorical interior operators were shown to be left adjoint neighbourhood operators (\cite{MR3196846}), a nice relationship between closure and neighbourhood operators in a category was lacking until the categorical topogenous structures (\cite{MR0286814}) and (\cite{iragi2016topogenous})  were recently introduced. Indeed the conglomerate of categorical topogenous structures is order isomorphic to the conglomerate all neighbourhood operators and contains both the conglomerates of all interior and all closure operators as reflective subcategories. From  this order isomorphism between neighbourhood operators and topogenous structures, one see that  all the concept studied using neighbourhood operators can be described using topogenous structures. However the use topogenous structures (particularly easier to work with) provides unified study of many concepts and notions that were separately studied for closure and interior operators. In the present paper wish to continue investigating the categorical topogenous structures and demonstrate their crucial role played in allowing a single setting study of several results obtained  separately for closure and interior operators on an abstract category.} 

{\indent  Classes of morphisms with respect to a closure $c$, interior $i$ and neighbourhood operators $\nu$ have been introduced and studied by a number of authors (see e.g \cite{MR1799498}). A morphism $f : X \longrightarrow Y$ is $c$-closed if it preserves the closure i.e $f(c(m)) = c(f(m))$ for any subobject $m$ of $X$.} An important question is then which morphism preserves a topogogenous order i.e $m\sqsubset n\Rightarrow f(m)\sqsubset f(n)$ for any subobject $m, n$ of $X$. Such morphism that we call $\sqsubset$-strict turns out to capture both the $c$-closed and  $i$-open morphisms. This leads to the definition of morphisms $f : X \longrightarrow Y$ which reflects a topogenous order: $f^{-1}(m)\sqsubset f^{-1}(n)\Rightarrow m\sqsubset n$ for any subobject $m, n$ of $Y$. We call it the $\sqsubset$-final morphism and show that it subsumes the $c$-final, the $i$-finial and $\nu$-finial. We also define the $\sqsubset$-co-strict and $\sqsubset$-initial. It is demonstrated that the $\sqsubset$-co-strict (resp. $\sqsubset$-initial) morphism provides a unified study of the $c$-open and  $i$-closed (resp. the $c$-initial and $i$-initial). While these four classes of morphisms are not closed under pulbacks, we show that each of them ascends along $\sqsubset$-initial morphisms and descends along $\sqsubset$-final morphisms. $\sqsubset$-Initial and $\sqsubset$-morphisms are shown to interact well with the $\sqsubset$-strict subobjects.

{\indent Topogenous structures induced by an $\mathcal{M}$-fibration are studied. They not only include D. Dikranjan and W. Tholen’s lifting of a closure as a particular case but also allow the lifting of an interior operator along this functor. We then define the continuity of a $\mathcal{C}$ -morphism with respect to two topogenous orders on $\mathcal{C}$ and use it to construct new topogenous orders from old. It is shown that for a pointed endofunctor $(F; \eta)$ of  $\mathcal{C}$  and a topogenous order on  $\mathcal{C}$, there is a coarsest topogenous order $ \sqsubset^{F,\eta}$  on $\mathcal{C}$  for which every $\eta_{X} : X \longrightarrow FX $ is $( \sqsubset^{F,\eta}, \sqsubset)$ -continuous and dually for a copointed endofunctor of $\mathcal{C}$, there is a finest topogenous order $ \sqsubset^{F,\eta}$  on $\mathcal{C}$ for which every $\varepsilon_{X} : GX \, \longrightarrow \ X $ is $(\sqsubset, \sqsubset^{G, \varepsilon})$ -continuous. In particular, for meet preserving topogenous orders, $ \sqsubset^{F, \eta} $ and $ \sqsubset^{G, \varepsilon} $ correspond to the closure operators obtained by Dikranjan and $ W $. Tholen in \cite{MR2122803} while for join preserving, $ \sqsubset^{F,\eta}$ and $ \sqsubset^{G,\varepsilon}$ allow us to construct the interior operators induced by $F$ and $G$ respectively.}

\section{Preliminaries} 
 {\indent  Our blanket reference for categorical concepts is \cite{MR2240597}. The basic facts on categorical closure operators used here can be found in \cite{MR2122803, dikranjan1987closure}. For the categorical topogenous structures, we use \cite{MR0286814, iragi2016topogenous}.
Throughout the paper, we consider a category $\mathcal{C}$ supplied with a proper $(\mathcal{E}, \mathcal{M})$-factorization system for morphisms. The category $\mathcal{C}$ is assumed to be 
$\mathcal{M}$-complete so that pullbacks of $\mathcal{M}$-morphisms along $\mathcal{C}$-morphisms and arbitrary $\mathcal{M}$-intersections of 
$\mathcal{M}$-morphisms exist and are again in $\mathcal{M}$. 
For any $X\in \mathcal{C}$, sub$X = \{m\in \mathcal{M}\;|\;\mbox{cod}(m) = X\}$. It is ordered as follows:
$n\leq m$ if and only if there exists $j$ such that $m\circ j = n$.
If $m\leq n$ and $n\leq m$ then they are isomorphic. We shall simply write $m = n$ in this case.
Sub$X$ is a (possibly large) complete lattice with greatest element $1_{X} : X\longrightarrow X$ and the least element $0_{X} : O_{X}\longrightarrow X$.}

 {\indent Any $\mathcal{C}$-morphism, $f:X\longrightarrow Y$  induces an image/pre-image adjunction 
 $f(m)\leq n$ if and only if $m\leq f^{-1}(n)$ for all $n\in \;$sub$Y$, $m\in \;$sub$X$ with $f(m)$ the $\mathcal{M}$-component of the $(\mathcal{E}, \mathcal{M})$-factorization of $f\circ m$ while $f^{-1}(n)$ is the pullback of $n$ along $f$. 
 We have from  the image/pre-image adjunction that $f(f^{-1}(n))\leq n$ (with $f(f^{-1}(n)) = n$ if $f\in \mathcal{E}$ and $\mathcal{E}$ is pullback stable along $\mathcal{M}$-morphisms) and $m\leq f^{-1}(f(m))$  (with $m = f^{-1}(f(m))$ if $f\in \mathcal{M}$)
for any $n\in \;$sub$Y$ and $m\in \;$sub$X$. }

 Applying adjointness repeatedly we obtain the lemma below.

\begin{Lemma}\label{L22}
 Let 
 
$$\xymatrix{  
       X'\ar[r]^{f'}  \ar[d]_{p'}  &  Y'\ar[d]^{p}\\
       X\ar[r]_{f} & Y
       }$$

be a commutative diagram. Then for any subobject $n\in \;$sub$Y'$,  $p'(f'^{-1}(n))\leq f^{-1}(p(n))$.
\end{Lemma}

\begin{defn}\label{D154}
 $A pointed \; endofunctor$ of $\mathcal{C}$ is a pair $(F,\eta)$ consisting of a functor $F : \mathcal{C}\longrightarrow \mathcal{C}$ and a natural transformation $\eta : 1_{\mathcal{C}}\longrightarrow F$.
 \end{defn}
For any $\mathcal{C}$-morphism $f : X\longrightarrow Y$, $(F,\eta)$ induces the commutative diagram below. 
$$\xymatrix{  
       X\ar[r]^{\eta_{X}}  \ar[d]_{f}  &  FX\ar[d]^{Ff}\\
       Y\ar[r]_{\eta_{Y}} & FY
       }$$ 
       
 If each $\eta_{X}\in \mathcal{F}$ where $\mathcal{F}$ is a class of  $\mathcal{C}$-morphisms, then $(F,\eta)$ is $\mathcal{F}$-pointed. 
 A copointed endofunctor of $\mathcal{C}$ is defined dually. 
\begin{defn} \label{D13}
	A topogenous order $\sqsubset$ on $\mathcal{C}$ is a family  $\{\sqsubset_{X}\; |\;X\in \mathcal{C}\}$ of relations, each
	$\sqsubset_{X}$ on sub$X$, such that:
	\begin{itemize}
		\item [$(T1)$] $m\sqsubset_{X} n\Rightarrow m\leq n$ for every $m, n\in \;$sub$X$,
		\item [$(T2)$] $m\leq n\sqsubset_{X} p\leq q\Rightarrow m\sqsubset_{X} q$ for every $m, n, p, q\in \;$sub$X$, and
		\item [$(T3)$] every morphism $f: X\longrightarrow Y$ in $\mathcal{C}$ is $\sqsubset$-continuous, $m\sqsubset_{Y} n\Rightarrow f^{-1}(m)\sqsubset_{X} f^{-1}(n)$ for every $m, n\in \;$sub$Y$.
	\end{itemize}
\end{defn}
Given two topogenous orders $\sqsubset$ and $\sqsubset'$ on $\mathcal{C}$, $\sqsubset\subseteq \sqsubset'$ if and only if
$m\sqsubset_{X} n\Rightarrow m\sqsubset'_{X} n$ for all $ m, n\in \;$sub$X$. The resulting ordered congolomerate of all topogenous orders on $\mathcal{C}$ is denoted by $TORD(\mathcal{C}, \mathcal{M}).$ A topogenous order $\sqsubset$ is said to be $\bigwedge$-preserving if $(\forall i\in I: m\sqsubset_{X} n_{i})\Rightarrow  m\sqsubset_{X} \bigwedge  n_{i}$. It is $\bigvee$-preserving if $(\forall i\in I: m_{i}\sqsubset_{X} n)\Rightarrow  \bigvee m_{i}\sqsubset_{X} n$. If $m\sqsubset_{X} n\Rightarrow \exists\;p\;|\;m\sqsubset_{X} p\sqsubset_{X} n$ for all $X\in \mathcal{C}$, then it is interpolative.
The ordered conglomerate of all $\bigwedge$-preserving and $\bigvee$-preserving topogenous orders is denoted by $\bigwedge-TORD(\mathcal{C}, \mathcal{M})$ and $\bigvee-TORD(\mathcal{C}, \mathcal{M})$ respectively. $\bigwedge-INTORD(\mathcal{C}, \mathcal{M})$(resp. $\bigwedge-INTORD(\mathcal{C}, \mathcal{M})$)  will denote the conglomerate of all interpolative meet preserving (resp. join preserving) topogenous orders.

\begin{defn}
	$A\;closure\;operator$ $c$ on $\mathcal{C}$ with respect to $\mathcal{M}$ is given by a
	family of maps \\$\{c_{X}$: sub$X\longrightarrow \;$sub$X\;|\; X\in \mathcal{C}\}$ such that:
	\begin{itemize}
		\item [$(C1)$]  $m\leq c_{X}(m)$ for all $m\in \;$sub$X$;
		\item [$(C2)$]  $m\leq n\Rightarrow c_{X}(m)\leq c_{X}(n)$ for all $m, n\in \;$sub$X$;
		\item [$(C3)$]   every morphism $f: X\longrightarrow Y$ is $c$-continuous: $f(c_{X}(m))\leq c_{Y}(f(m))$ for all
		$m\in \;$sub$X$.
	\end{itemize}
\end{defn}

We denote by $CL(\mathcal{C}, \mathcal{M})$ the conglomerate of all closure operators on $\mathcal{C}$ with respect to $\mathcal{M}$ ordered as follows: 
$c\leq c'$ if $c_{X}(m)\leq c'_{X}(m)$ for all $m\in \;$sub$X$ and $X\in \mathcal{C}$.  A closure operator $c$ on $\mathcal{C}$ is idempotent if $c_{X}(c_{X}(m)) = c_{X}(m)$ for all $m\in \;$sub$X$ and $X\in \mathcal{C}$. $ICL(\mathcal{C}, \mathcal{M})$ will denote the conglomerate of all idempotent closure operators on $\mathcal{C}$.

\begin{defn}
	$An\;interior\;operator$ $i$ on $\mathcal{C}$ with respect to $\mathcal{M}$ is given by a
	family of maps \\$\{i_{X}$: sub$X\longrightarrow \;$sub$X\;|\; X\in \mathcal{C}\}$ such that:
	\begin{itemize}
		\item [$(I1)$]  $i_{X}(m)\leq m$ for all $m\in \;$sub$X$ and $X\in \mathcal{C}$;
		\item [$(I2)$]  $m\leq n\Rightarrow i_{X}(m)\leq i_{X}(n)$ for all $m, n\in \;$sub$X$, $X\in \mathcal{C}$;
		\item [$(I3)$]   every morphism $f: X\longrightarrow Y$ is $i$-continuous: $f^{-1}(i_{Y}(n))\leq i_{X}(f^{-1}(n))$ for all $n\in \;$sub$Y$.
	\end{itemize}
\end{defn}
$INT(\mathcal{C}, \mathcal{M})$ will denote the conglomerate of all interior operators on $\mathcal{C}$ with respect to $\mathcal{M}$. It is ordered as follows: $i\leq  i'$ if ¨$i_{X}(m)\leq i'_{X}(m)$ for all $m\in \;$sub$X$, $X\in \mathcal{C}$.  A subobject of an object $X\in \mathcal{C}$ is said to be $i$-open if $i_{X}(m) = m$. An interior operator $i$ is said to be idempotent provided that $i_{X}(m)$ is $i$-open for every $m\in \;$sub$X$ and $X\in \mathcal{C}$. We write $IINT(\mathcal{C}, \mathcal{M})$ the ordered conglomerate of all idempotent interior operators on $\mathcal{C}$.
\begin{defn}
	A neigbhourhood operator $\nu$ on $\mathcal{C}$ with respect to $\mathcal{M}$ is given by a
	family of maps $\{c_{X}$: sub$X\longrightarrow \;P($sub$X)\;|\; X\in \mathcal{C}\}$ such that:
	\begin{itemize}
		\item [$(N1)$]  $n\in \nu_{X}(m)\Rightarrow m\leq n$ for all $m, n\in \;$sub$X$;
		\item [$(N2)$]  $m\leq n\Rightarrow \nu_{X}(n)\leq \nu_{X}(m)$ for all $m, n\in \;$sub$X$;
		\item [$(N3)$]  $p\in \nu_{X}(m)$ and $p\leq q$ then $q\in \nu_{X}(m)$ for all $p, q, m\in \;$sub$X$;
		\item [$(N4)$]  every morphism $f: X\longrightarrow Y$ is $\nu$-continuous: $n\in \nu_{Y}(f(m))\Rightarrow f^{-1}(n)\in \nu_{X}(m)$ for all
		$m\in \;$sub$X$ and $n\in \;$sub$Y$.
	\end{itemize}
\end{defn}
The congolomerate of all neighbourhood operators on $\mathcal{C}$ with respect to $\mathcal{M}$ is denoted by
 $NBH(\mathcal{C}, \mathcal{M})$.
 
 {\indent$\;\;\;\;$ While associated  closure and interior operators provide an equivalent description a topology. This natural correspondance between these two notions no longer hold in abstract categorical settings since the subobject lattices are not necessarly boolean algebras. The categorical topogenous structure (particularly easier to work with) that is defined above allows to obtain a nice relationship between closure and interior operators on a category in the sense that many concepts and definitions that have been studied separately  for categorical closure and interior operators  can be shown to be exactly the same using  topogenous orders (see for example the classes of morphisms definied and studied in this paper). The next result that we recall from ((\cite{MR0286814}, \cite{iragi2016topogenous})) exibits the clear relationship between closure, interior and topogenous order in a category.}
 
\begin{Proposition}\label{P20}
	
$\;\;\;\;\;\;\;\;\;\;\;\;\;\;\;\;\;\;\;\;\;\;\;$
\begin{itemize}
 \item [$(1)$] $TORD(\mathcal{C}, \mathcal{M})$  and  $NBH(\mathcal{C}, \mathcal{M})$ are order isomorphic with the inverse assignments $\sqsubset\longrightarrow \nu^{\sqsubset}$ and $\nu\longrightarrow \sqsubset^{\nu}$  given by $\nu^{\sqsubset}_{X}(m) = \{n\; |\; m\sqsubset_{X} n\}\;\mbox{and}\;\; m\sqsubset^{\nu}_{X} n\Leftrightarrow n\in \nu_{X}(m)\;\mbox{for all}\;X\in \mathcal{C}.$
	
 \item [$(2)$] $\bigvee-TORD(\mathcal{C}, \mathcal{M})$ is order isomorphic to $INT(\mathcal{C},\mathcal{M}).$
 The inverse assignments  of each other are given by $i^{\sqsubset}_{X}(m) = \bigvee\{p\; |\; p\sqsubset_{X} m\}\;\mbox{and}\;\; m\sqsubset^{i}_{X} n\Leftrightarrow m\leq i_{X}(n)\;\mbox{for all}\;X\in \mathcal{C}.$
\item [$(3)$]$\bigwedge-TORD(\mathcal{C}, \mathcal{M})$ is order isomorphic to $CL(\mathcal{C},\mathcal{M})$. The inverse assignments  of each other are given by $c^{\sqsubset}_{X}(m) = \bigwedge\{p\; |\; m\sqsubset_{X} p\}\;\mbox{and}\;\; m\sqsubset^{c}_{X} n\Leftrightarrow c_{X}(m)\leq n\;\mbox{for all}\;X\in \mathcal{C}.$
 	\end{itemize} 	 
 \end{Proposition}

\section{Families of morphisms}

\begin{Proposition}
	Assume that $f^{-1}$ has a right adjoint for any $f \in \mathcal{C}$. Let $f : X \longrightarrow Y$ be a
	$\mathcal{C}$-morphism and $\sqsubset \in TORD(\mathcal{C}, \mathcal{M})$. The following are equivalent to the $\sqsubset$-continuity of $f$.
	\begin{enumerate}
		\item [(1)] $m \sqsubset_Y n \Rightarrow f^{-1}(m) \sqsubset_X f^{-1}(n), m \in subY \ and \ n \in subX$;
		\item [(2)] $m \sqsubset_Y f_{*}(n) n \Rightarrow f^{-1}(m) \sqsubset_X n, m \in subY \ and \ n \in subX$;
		\item[(3)] $m \sqsubset_Y f_{*}(n) n \Rightarrow m\sqsubset_X n, m,n \in subY \ and \ n \in subX$;
	\end{enumerate}
\end{Proposition}
\begin{proof}
If $(T 3)$ holds, then $m \sqsubset_X n \Rightarrow  f(f^{-1}(m)) \leq m \sqsubset_X n \Rightarrow f(f^{-1}(m)) \sqsubset_X n \Rightarrow 
f^{-1}(m) \sqsubset_X f^{-1}(n)$. Assume (1) holds, then $m \sqsubset_Y f_{*}(n) \Rightarrow m \leq f^{-1}(f(m)) \sqsubset_X n \Rightarrow m \sqsubset_X n$.
If (2) holds then $f(m) \sqsubset_Y f_{*}(n) \Rightarrow m \leq f^{-1}(f(m)) \sqsubset_X n \Rightarrow m \sqsubset_X n$.
If (3) holds, then $f(m) \sqsubset_Y p \Rightarrow f(m) \sqsubset_X p \leq f_{*}(f^{-1}(p)) \Rightarrow f(m) \sqsubset_Y f_{*}(f^{-1}(p)) \Rightarrow
m \sqsubset_X f^{-1}(p)$. 
\end{proof}
The point in the proposition above that each equivalent description of $(T3)$ fullfils one implication leads to the kind of morphisms that satisfy the other implication.

\begin{defn}
	Let $f^{-1}$ have a right adjoint for any $f \in \mathcal{C}$ and $\sqsubset \in TORD$. We say that a $\mathcal{C}$-morphism $f : X \longrightarrow Y$ is
	\begin{enumerate}
		\item [(1)] [(13)] $\sqsubset$-strict if $f(m) \sqsubset_Y p \Leftrightarrow m \sqsubset_{X} f^{-1}(p) \ for \ all \ m \in subX \ and \ p \in subY$
		\item [(2)] $\sqsubset-final$ if $m \sqsubset_Y n \Leftrightarrow f^{-1}(m) \sqsubset_{X} f^{-1}(n)$ for all $n, m \in subY$ ;
		\item [(3)] $\sqsubset$-co-strict if $m \sqsubset_Y f_{*}(n) , f^{-1}(m) \sqsubset_X n$ for all $m \in subY \ and \ n \in subX$,
		\item [(4)] $\sqsubset$-initial if $f(m) \sqsubset_Y f_{*}(n) \Leftrightarrow m \sqsubset_X n $ for all $m, n \in subX$.
	\end{enumerate}
\end{defn}
It follows immediately from Proposition \ref{P20}$(1)$ that a morphisme $f : X \longrightarrow Y$ is $\sqsubset$-strict (resp. $\sqsubset$-co-strict, $\sqsubset$-final, $\sqsubset$-initial) if and only if $\nu^{\sqsubset}_{Y} (f(m)) = f(\nu^{\sqsubset}_{X}(m))$ (resp. $\nu^{\sqsubset}_{X} (f^{-1}(n)) = f^{-1}(\nu^{\sqsubset}_{Y}(n))$, $f(\nu^{\sqsubset}_{X} (f(n))) = \nu^{\sqsubset}_{Y}(n)$, $f^{-1}(\nu^{\sqsubset}_{Y} (f(m))) = \nu^{\sqsubset}_{X}(m)$). These morphisms were defined and in (\cite{razafindrakoto2013neighbourhood}) with respect to a neighbourhood operator. Recall from (\cite{iragi2016topogenous}) that a subobject m of an object $X \in \mathcal{C}$ is $\sqsubset$-strict if $m \sqsubset_X m$.

\begin{Proposition}
	Let $f : X \longrightarrow Y$ be a $\mathcal{C}$-morphism.
	\begin{enumerate}
		\item[(1)] If $f$ is $\sqsubset$-final then a subobject $n$ of $Y$ is $\sqsubset$-strict $iff f^{-1}(n)$ is $\sqsubset$-strict in $X$.
		\item[(2)] If $\sqsubset \in INTORD(\mathcal{C}, \mathcal{M})$ and $f$ is $\sqsubset$-initial then for every $\sqsubset$-strict subobject $m$ of $X$,
		there is $p \in subY$ such that $m = f^{-1}(p)$.
	\end{enumerate}	
\end{Proposition}

\noindent \textit{Proof.}
(1) is clear. For (2), assume $\sqsubset \in INTORD(\mathcal{C}, \mathcal{M})$, $f$ is $\sqsubset$-initial and $m \in subX$ is $\sqsubset$-strict. Then $m \sqsubset_X m \Rightarrow f(m) \sqsubset_X f_{*}(m) \Rightarrow \exists \ p \in subX$ such that $f(m) \sqsubset_X p \sqsubseteq_X f_{*}(m) \Rightarrow m \sqsubset_{X} f^{-1}(p) \sqsubset_X f^{-1}(f_{*}(m)) \leq m \Rightarrow m \leq f^{-1}(p) \leq m \Rightarrow m = f^{-1}(p). \qquad \qquad \square$ \\ 

\par The following easy proposition to prove summarizes the behaviour of our classes of morphisms introduced in Definition 3.1.

\begin{Proposition}
	\begin{enumerate}
		\item[(1)] Each of the classes is closed under composition and contains all
		isomorphisms of $\mathcal{C}$.
		\item [(2)] $\sqsubset$-initial morphisms are left-cancelable, while $\sqsubset$-co-strict, $\sqsubset$-strict and <-final morphisms are left cancellable with respect to $\mathcal{M}$.
		\item [(3)] $\sqsubset$-final morphisms are right cancelable, while $\sqsubset$-initial, $\sqsubset$-co-strict and $\sqsubset$-strict morphisms are right cancellable with respect to $\varepsilon$ provided $\varepsilon$ is pullback stable.
	\end{enumerate}
\end{Proposition}
The following is an observation concerning the relationship between the types of morphisms.
\begin{Proposition}
	\item [(1)] Every $\sqsubset$-co-strict morphism in $\mathcal{M}$ is $\sqsubset$-initial.
	\item [(2)] Every $\sqsubset$-initial morphism in $\varepsilon$ is $\sqsubset$-co-strict provided $\varepsilon$ is pullback stable.
	\item [(3)] Any $\sqsubset$-strict morphism in $\mathcal{M}$ is $\sqsubset$-initial.
	\item [(4)] Every $\sqsubset$-strict in $\varepsilon$ is $\sqsubset$-final provided $\varepsilon$ is pullback stable.
	\item [(5)] Any $\sqsubset$-final morphism in $\mathcal{M}$ is $\sqsubset$-strict.
	\item [(6)] Every $\sqsubset$-co-strict morphism in $\varepsilon$ is $\sqsubset$-final provided $\varepsilon$ is pullback stable.
	\item [(7)] If $g \circ f = 1$ in $\mathcal{C}$ then $f$ is a $\sqsubset$-initial morphism and g is a $\sqsubset$-final morphism in $\varepsilon$.
	
\end{Proposition}
We are interested in the pullback behaviour of the morphisms. We show that each of the classes ascends along $\sqsubset$-initial morphisms and descends along $\sqsubset$-final morphisms
\begin{Proposition}
Let  
$$\xymatrix{  
	X'\ar[r]^{f'}  \ar[d]_{p'}  &  Y'\ar[d]^{p}\\
	X\ar[r]_{f} & Y
}$$ 
	be a pullback diagram satisfying the Beck Chevalley’s Property. Then the following statements
	are true.
	\begin{enumerate}
		\item [(i)] If $p^{'}$ is $\sqsubset$-initial, then $f^{'}$ is $\sqsubset$-initial (resp. $\sqsubset$-strict, $\sqsubset$-co-strict, $\sqsubset$-final) provided f
		is $\sqsubset$-initial (resp. $\sqsubset$-strict, $\sqsubset$-co-strict, $\sqsubset$-final).
		\item[(ii)] If $p$ is $\sqsubset$-final, then $f$ is $\sqsubset$-final (resp. $\sqsubset$-strict, $\sqsubset$-co-strict, $\sqsubset$-initial) provided $f^{'}$ is
		$\sqsubset$-final (resp. $\sqsubset$-strict, $\sqsubset$-co-strict, $\sqsubset$-initial).
	\end{enumerate}
\end{Proposition}

\begin{proof}
A sample of the calculations is given below.
\begin{enumerate}
	\item[(i)] Let $p^{'}$ be $\sqsubset$-initial and assume that f is $\sqsubset$-strict. Then by using the commutativity of the diagram, $\sqsubset$-continuity of $p$ and BCP, we have that  
	$m \sqsubset n \Rightarrow p^{'}(m)\sqsubset_{X}p^{'}_{*}(n)\Rightarrow f(p^{'}(m))\sqsubset_Y f(p^{'}_{*}(n))\Rightarrow p(f^{'}(m))\sqsubset_Y f(p^{'}_{*}(n))\Rightarrow f^{'}(m)\sqsubset_Q p^{-1}(f(p^{'}_{*}(n)))\Rightarrow f^{'}(m)\sqsubset_Q f^{'}(p^{'-1}(p^{'}_{*}(n))) \leq f^{'}(n) \Rightarrow f^{'}(m) \sqsubset_Q f^{'}(n)$.
	
	\item [(ii)] Let $p$ be final and $f^{'} \sqsubset$-initial. Then by using the $\sqsubset$-continuity of $p^{'}$ and BCP, we have that
                $m \sqsubset_X n \Rightarrow p^{'-1}(m)\sqsubset_{P}p^{'-1}(m) \Rightarrow f^{'}(p^{'-1}(m))\sqsubset_Q f^{'}_{*}(p^{'-1}(n))\Rightarrow 
		p^{-1}(f(m))\sqsubset_Q p^{-1}(f_{*}(n))\Rightarrow f(m) \sqsubset_Y f_{*}(n)$.	
	\end{enumerate}
\end{proof}
Our next proposition shows that when $\sqsubset \in \bigwedge-TORD(\mathcal{C}, \mathcal{M})$, then the $\sqsubset$-initial morphism (resp. $\sqsubset$-co-strict, $\sqsubset$-final, $\sqsubset$-strict ) correspond the $c$-initial (resp. $c$-open, $c$-final,
$c$-closed ) morphisms that have been studied in the literature (see e.g \cite{MR2122803, MR1799498}).
\begin{Proposition}
	Let $\sqsubset \in \bigwedge-TORD (\mathcal{C}, \mathcal{M})$ and$ c \in CL (\mathcal{C}, \mathcal{M})$. Let for any morphism $
	f \in \mathcal{C}$, the inverse image $f^{-1}$ commutes with the join of subobjects. Then $f : X \longrightarrow Y$
	is $\sqsubset $-initial (resp. $\sqsubset $-co-strict, $\sqsubset $-final, $\sqsubset $-strict ) if and only if it is $c$-initial (resp $c$-open,
	$c$-final, $c$-closed ). 
\end{Proposition}

\begin{proof} 
\begin{enumerate}
	\item [(1)] Let f be $c$-initial. Then $f(m) \sqsubset^{c}_{Y} f_{*}(n) \Leftrightarrow c^{\sqsubset}_{Y}(f(m)) \leq f_{*}(n) \Leftrightarrow f^{-1}(c^{\sqsubset}_{Y}(f(m))) \leq n \Leftrightarrow c^{\sqsubset}_{X}(m) \leq n \Leftrightarrow m \sqsubset^{c}_{X} n$. Conversely if $f$ is $\sqsubset$-initial, then $f^{-1}(c^{\sqsubset}_{Y}(f(m)) \leq n \Leftrightarrow
	c^{\sqsubset}_{Y}(f(m)) \leq f_{*}(n) \Leftrightarrow f(m) \sqsubset_{Y} f_{*}(n) \Leftrightarrow m \sqsubset_{X} n \Leftrightarrow c^{\sqsubset}_{Y}(m) \leq n$.
	\item [(2)] Assume $f$ is $c$-open, then $m \sqsubset^{c}_{Y} f_{*}(n) \Leftrightarrow \sqsubset^{c}_{Y} (m) \leq f_{*}(n) \Leftrightarrow f^{-1}(c^{\sqsubset}_{Y} (m) \leq n \Leftrightarrow
	c^{\sqsubset}_{X}(f^{-1}(m)) \leq n \Leftrightarrow f^{-1}(m) \sqsubset^{c}_{X} n$. Conversely if $\sqsubset$-co-strict then $c^{\sqsubset}(f^{-1}(m)) \leq n \Leftrightarrow
	f^{-1}(m) \sqsubset_{X} n \Leftrightarrow m \sqsubset_{Y} (n) \Leftrightarrow c^{\sqsubset}_{Y} (m) \ge f_{*}(n) \Leftrightarrow f^{-1}(c^{\sqsubset}_{Y} (m)) \leq n$.
	
	\item [(3)]  If f is $c$-final, then $f^{-1}(m) \sqsubset^{c}_{X} f^{-1}(n) \Leftrightarrow c^{\sqsubset}(f^{-1}(m)) \leq f^{-1}(n) \Leftrightarrow f(c^{\sqsubset}(f^{-1}(m)) \leq
	n \Leftrightarrow c_Y (m) \leq n \Leftrightarrow m \sqsubset^{c}_{Y} n$. On the other hand if $f$ is $\sqsubset$-final, then $f(c^{\sqsubset}(f^{-1}(m)) \leq
	n \Leftrightarrow c^{\sqsubset}_{X}(f^{-1}(m)) \leq f^{-1}(n) \Leftrightarrow f^{-1}(m) \sqsubset_{X} f^{-1}(n) \Leftrightarrow m \sqsubset_{Y} n \Leftrightarrow c^{\sqsubset}_{Y} (m) \leq n$.
	
	\item [(4)] See (\cite{MR0286814}, Propostion $3.5$).
\end{enumerate}
\end{proof}
Open morphisms with respect to an interior operator were studied in (\cite{MR3351083}). Assuming the pre-image commutes with the join of subobjects, $i$-initial and $i$-final morphisms were introduced in \cite{razafindrakoto2013neighbourhood} . Recently in \cite{assfaw2019interior}, the $i$-closed morphism has been defined and a systematic study of the four classes of morphisms with respect to an interior operator is provided. In the next proposition we prove that if $\sqsubset \in \bigvee -TORD(\mathcal{C}, \mathcal{M})$, then the $\sqsubset$-strict morphism (resp.  $\sqsubset$-co-strict,  $\sqsubset$-initial,  $\sqsubset$-final) correspond to the $i$-open (resp. $i$-closed, $i$-initial and $i$-final) morphisms. 
\begin{Proposition}
	Let $\sqsubset \in \bigvee-TORD(\mathcal{C},\mathcal{M})$, and let for any morphism $f \in \mathcal{C}$, the preimage $f^{-1}$ commute with the join of subobjects. Then $f : X \longrightarrow Y$ is $\sqsubset$-initial (resp.
	$\sqsubset$-co-strict, $\sqsubset$-final, $\sqsubset$-strict) if and only if it is $i$-initial (resp. $i$-closed, $i$-final, $i$-open ).
\end{Proposition}
Apart from the four classes of morphisms studied above, a weaker notion of <-final
morphism will be useful.
\begin{defn}
	Let $\sqsubset \in TORD$. A $\mathcal{C}$-morphism $f : X \longrightarrow Y$ is said to be weakly $\sqsubset$-final
	if for any $m, n \in subY$ such that $m \leq n, m \sqsubset_{Y} n \Rightarrow f^{-1}(m) \sqsubset_{X} f^{-1}(n)$.
	\end{defn}
We note that if $f \in \varepsilon$, then $f$ is weakly $\sqsubset$-final if and only if it is $\sqsubset$-final.
\begin{Proposition}
	Let $\sqsubset \in TORD$ and $f : X \longrightarrow Y$ be a $\mathcal{C}$-morphism.
	\begin{enumerate}
		\item [(1)] If $\sqsubset \in \bigwedge -TORD$, then $f$ is weakly $\sqsubset$-final if and only if $c^{\sqsubset}_{Y} (m) = m \vee f(c^{\sqsubset}_{X}(f^{-1}(m)))$.
		\item [(2)] If $ \sqsubset \in \bigvee -TORD$, then $f$ is weakly $\sqsubset$-final if and only if $i^{\sqsubset}_{Y} (m) = m\wedge f_{*}(i^{\sqsubset}_{X} (f^{-1}(m)))$.	
	\end{enumerate}
	\end{Proposition}

\section{Lifting a Topogenous order along an $\mathcal{M}$-fibration} 
Recall from (\cite{MR2122803}) that for an $\mathcal{M}$-fibration $F : \mathcal{A} \longrightarrow \mathcal{C}, (\varepsilon_F ,\mathcal{M}_F ) \mbox{ where } \varepsilon_F = F^{-1} \varepsilon = \{e \in
\mathcal{C} | F_e \in \varepsilon \} \mbox{ and } \mathcal{M}_F = F^{-1} \mathcal{M} \bigcap \textit{IniF}$, with \textit{IniF} the class of \textit{F}-initial morphisms, is a factorization system in $\mathcal{A}$ and $\mathcal{M}$-subobject properties in $\mathcal{C}$ are inherited by $\mathcal{M}_F$-subobjects in $\mathcal{A}$. In particular, the $\mathcal{M}_F$-images and $\mathcal{M}_F$-pre-images are obtained by y initially lifting
$\mathcal{M}$-images and $\mathcal{M}$-inverse images. Consequently $Ff^{-1}(n) = (Ff)^{-1}(Fn)$ and $(Ff)(Fm) = Ff(m)$ far any $f \in \mathcal{A}$ and suitable subobjects $n$ and $m$.  Let $F : \mathcal{A} \longrightarrow \mathcal{C}$ be a faithful $\mathcal{M}$-fibration. For any $X \in \mathcal{A}$, sub$X$ and sub$FX$ are order equivalent with the inverse assignments, $\gamma_X : subX \longrightarrow subFX \mbox{ and } \delta_X : subFX \longrightarrow subX$, given by $\gamma_X(m) = Fm \mbox{ and } \delta_X(n) = p \mbox{ with } Fp = n \mbox{ and } p \in IniF$.  For any $f : X \longrightarrow Y \in \mathcal{A}$ and suitable  subobjects $n, m, n' \mbox{ and }m'$.  $\gamma_Y(f(m)) = (Ff)(\gamma_X(m))$, $f(\delta_X(n)) = \delta_Y(Ff)(n)$, $f^{-1}(\delta_Y(m')) = \delta_X((Ff)^{-1}(m'))$ and $\gamma_X(f^{-1}(n')) = (Ff)^{-1}(\gamma_Y(n'))$.

 \begin{Proposition}
 	Let $F : \mathcal{A} \longrightarrow \mathcal{C}$ be a faithful $\mathcal{M}$-fibration and $\sqsubset$ be a topogenous order on $\mathcal{C}$ with respect to $\mathcal{M}$. Define $\sqsubset^F \mbox{ by } m \sqsubset^F_X n \Leftrightarrow Fm \sqsubset_{FX} \gamma_X (n) = Fn.$ Then $\sqsubset^F$ is a topogenous order on $\mathcal{A}$ with respect to $\mathcal{M}_F.$ Moreover $\sqsubset^F$ is interpolative and satisfies (T4) provided $\sqsubset$ has the same properties.
 \end{Proposition}
  \begin{proof}
   (T1) and (T2) are clear. For (T3), Let $f : X \longrightarrow Y$ be a $\mathcal{A}$-morphism and $f(m) \sqsubset^F_X n. \mbox{ Then } Ff(m) \sqsubset_{FY} \gamma_Y (n) \Rightarrow (Ff)(Fm) \sqsubseteq_{FY}(Ff)^{-1}(\gamma_Y(n)) = \gamma_X (f^{-1}(n)) \Leftrightarrow m \sqsubset^F_X f^{-1}(n).$
    If $\sqsubset^F$ is interpolative, then $m \sqsubset^F_X n \Leftrightarrow Fm \sqsubset_{FX} \gamma_X (n) \Rightarrow \exists p \in subFX | Fm \sqsubset_{FX} p \sqsubset_{FX} \gamma_X (n) \Rightarrow Fm \sqsubset_{FX} \gamma_X(\delta_X(p)) \mbox{ and } F(\gamma_X (p)) \sqsubset_{FX} \gamma_X (n), \mbox{ since } F\delta_X (p) = p \mbox{ and } m' \sqsubset^F_X n'. \mbox{ This implies that } Fm \sqsubset_{FX} \gamma_X (n) \mbox{ and } Fm' \sqsubset_{FX} (\gamma_X (n)). \\ \mbox{ Thus } Fm \bigwedge Fm' \sqsubset_{FX} \gamma_X (n) \bigwedge \gamma_X (n') \Rightarrow F(m \bigwedge m') \sqsubset_{FX} \gamma_X (n \bigwedge n') \Leftrightarrow m \bigwedge m' \sqsubset^{FX}_X n \bigwedge n'.$
    \end{proof}    
 In the light of Proposition 2.3, we can prove the following.
 
 \begin{Proposition}
 \begin{enumerate}
 	\item $If \sqsubset \in \bigwedge - TORD \mbox{, then } m \sqsubset^F_X n \Leftrightarrow \delta_X (\mathcal{C}^\sqsubset_{FX}(Fm)) \leq n$
    \item  $If \sqsubset \in \bigvee - TORD \mbox{, then } m \sqsubset^F_X n \Leftrightarrow m \leq \delta_X (i^\sqsubset_{FX}(\gamma_X (n))).$
   \end{enumerate}   
 \end{Proposition}
 \begin{proof}
 	 \begin{enumerate}
 	 	\item $If \sqsubset \in \bigwedge - TORD \mbox{then} m \sqsubset^F_X n \Leftrightarrow Fm \sqsubset_{FX} \gamma_X (n) \Leftrightarrow c_{FX}(Fm) \leq \gamma_X \Leftrightarrow \delta_X c^\sqsubset_{FX}(Fm)) \leq \delta_X (\gamma_X (n)) \Leftrightarrow \delta_X (c^\sqsubset_{FX}(Fm)) \leq n.$
   	    \item  $If \sqsubset \in \bigvee - TORD \mbox{, then } m \sqsubset^F_X n \Leftrightarrow Fm \sqsubset_{FX} \gamma_X (n) \Leftrightarrow Fm \leq i^\sqsubset_{FX}(\gamma_X (n)) \Leftrightarrow \delta_X (Fm) \leq \delta_X (i^\sqsubset_{FX}(\gamma_X (n))) \Leftrightarrow m \leq delta_X (i^\sqsubset_{FX}(\gamma_X (n)))$.
   \end{enumerate}
\end{proof}
  \section{Topogenous orders induced by (co)pointed endofunctors}	
  \begin{defn}
 Let $\sqsubset, \sqsubset ' \in TORD(\mathcal{C}, \mathcal{M}).$ A $\mathcal{C}$-morphism $f : X \longrightarrow Y$  is $(\sqsubset, \sqsubset ')$- continuous if $f(m) \sqsubset^{'}_Y n \Rightarrow m \sqsubset_X f^{-1} (p) \sqsubset_X f^{-1}(n) \mbox{ for all } n, p \in $sub$Y \mbox{and} m \in $sub$X.$
  \end{defn}
  {\indent It is clear from the definition that every $\mathcal{C}-$ morphism $ f : X \longrightarrow Y$ is $(\sqsubset, \sqsubset)$-continuous and  $(\sqsubset ', \sqsubset')$-continuous, it is $(\sqsubset, \sqsubset')$-continuous if $\sqsubset' \subseteq \sqsubset$. If $\sqsubset, \sqsubset' \in \bigwedge-TORD(\mathcal{C}, \mathcal{M})$, then $f : X \longrightarrow Y$ is $ (\sqsubset, \sqsubset')$-continuous if and only if $f^{-1}(c^\sqsubset_X (m)) \leq c^\sqsubset_X (m)$ (will be referred to as ( $c^\sqsubset, c^{\sqsubset^{'}}$)-continuous). If $\sqsubset, \sqsubset' \in \bigvee-TORD(\mathcal{C}, \mathcal{M})$, then $f : X \longrightarrow Y$ is $(\sqsubset, \sqsubset')$-continuous if and only if $ f^{-1}(i^\sqsubset_Y (n)) \leq i^\sqsubset_X (f^{-1}(n))$(will be referred to as ($i^\sqsubset, i^{\sqsubset^{'}}$)-continuous). }
  	
{\indent For a pointed endofunctor ($F; \eta$) of $\mathcal{C}$ and a topogenous order $\sqsubset$ on $\mathcal{C}$, we wish to construct
the coarsest topogenous order $\sqsubset '$ on $\mathcal{C}$ for which every morphism in $\mathcal{F} = \{\eta_Y : X \in \mathcal{C}\}$ is
($\sqsubset^{'}$, $\sqsubset$)-continuous and the dual case. This method was developped for categorical closure operators in \cite{MR2122803} and it is used in our paper (\cite{MR0393425}) in the case of categorical quasi-uniform structures. The interested reader is referred to (\cite{MR0396825, iragi2019quasi, holgate2019quasi, iragi2023transitive}) for categorical quasi-uniform structures. In \cite{iragi2023topogenous}, topogenous orders have been recently introduced on faithful and amnestic functors. 
 \begin{Theorem}
 	Let ($F; \eta$) be an $\varepsilon$-pointed endofunctor of $\mathcal{C}$ and $\sqsubset$ be a topogenous order on $\mathcal{C}$. Then for all $m, n \in subX, m \sqsubset^{F, \eta}_X n \Leftrightarrow \eta_X (m) \sqsubset_{FX} p$ and $\eta^{-1}_X (p) \leq n$ is a topogenous
 	order on $\mathcal{C}$. It is the least topogenous order for which every $\eta_X : X \longrightarrow FX$ is ($\sqsubset^{F, \eta}, \sqsubset$)-
 	continuous. Moreover, $\sqsubset^{F, \eta}$ is interpolative provided $\sqsubset$ interpolates.
 \end{Theorem}
 \begin{proof} 
 	$(T1)$ and $(T2)$ is easily seen to be satisfied. To check $(T3)$, let $X \longrightarrow Y$ be a $\mathcal{C}$ mor-
 phism and $f(m) \sqsubset^{F, \eta}_Y n$. Then there is $p \in subFY$ such that $\eta_Y (m) \sqsubset_{FY} p$ and $\eta^{-1}_Y (p) \leq n$.
 By Definition 2.1, $Ff \circ \eta_X = \eta_Y \circ f$. So $(Ff)(\eta_X (m)) \sqsubset_FY p$ and $\eta^{-1}_Y (p) \leq n \Rightarrow (Ff)(\eta_X (m)) \sqsubset_FY p$ and $f^{-1}(\eta^{-1}_Y (p)) \leq f^{-1}(n) \Rightarrow \eta_X (m)\sqsubset_X (Ff)^{-1}(p)$ and $\eta^{-1}_X ((Ff)^{-1}(p)) \leq f^{-1}(n) \Rightarrow \eta_X (m) \sqsubset_{FX} l$ and $\eta_X^{-1}(l) \leq f^{-1}(n)(l = g^{-1}(p)) \Leftrightarrow m \sqsubset^{F, \eta}_X f^{-1}(n)$. Since $\eta_X (m) \sqsubset_{FX} n \Rightarrow \eta_X (m) \sqsubset_{FX} n \leq \eta_X (\eta^{-1}_X (n)) \Rightarrow \eta_X (m) \sqsubset_{FX} \eta_X (\eta_X^{-1}(n)) \Leftrightarrow m \sqsubset^{F, \eta}_X \eta^{-1}_X (n), F$ is $(\sqsubset^{F,\eta}, \sqsubset)$-continuous. If $\sqsubset^{'}$ is another topogenous order on $\mathcal{C}$ such that $F$ is $(\sqsubset^{'}, \sqsubset)$-continuous, then $m \sqsubset^{F, \eta}_X n \Leftrightarrow \eta_X (m) \sqsubset_{FX} p$ and $\eta_X^{-1}(p) \leq n \Rightarrow m \sqsubset_X \eta_X^{-1}(p) \leq n \Rightarrow m \sqsubset^{'}_X n$. It is also not hard to see that $\sqsubset^{F, \eta}$ interpolates provided $\sqsubset$ has the same property.
  \end{proof} 
 Since a reflector can be viewed as pointed endofunctor, one obtains the following proposition that will turn out to be useful in the examples.
\begin{Proposition}
	Let $\mathcal{A}$ be an $\varepsilon$-reflective subcategory of $\mathcal{C}$ and $\sqsubset$ be a topogenous order on
	$\mathcal{A}$. Then for all $X \in \mathcal{C}$ and $m, n \in subX, m \sqsubset^\mathcal{A} n \Leftrightarrow \eta_X (m) \sqsubset_{FX} p$ and $\eta_X^{-1} (p) \leq n$ is a topogenous order on $\mathcal{C}$. It is the least one for which every reflection morphism $\eta_X : X \longrightarrow	FX$ is $(\sqsubset^{F, \eta}, \sqsubset)$-continuous. Moreover, $\sqsubset^{F, \eta}$ is interpolative provided $\sqsubset$ interpolates.
\end{Proposition}

The next proposition is obtained from Proposition $2.3(b)$.

\begin{Proposition}
	Let $(F; \eta)$ be a pointed endofunctor of $\mathcal{C}$ and $\sqsubset \in \bigwedge-TORD.$ \\ Then $c^{\sqsubset^{F, \eta}}(m) = \eta^{-1}(c^{\sqsubset}_{FX}(\eta_X (m)))$ is the largest closure operator on $\mathcal{C}$ for which every $\eta_X : X \longrightarrow FX$ is $(c^{\sqsubset^{F, \eta}}, c^{\sqsubset})$-continuous. Moreover, if $\sqsubset$ is interpolative, then $c^{\sqsubset^{F, \eta}}$ is idempotent.
\end{Proposition}
\begin{proof} We easily check $(C1)$ and $(C2)$ is clear. For $(C3)$, Let $f : X \longrightarrow Y$ be a $\mathcal{C}$-morphism and $m \in subX$. Then from $\mathcal{C}^{\sqsubset}$-continuity, Lemma 2.1 and Definition 2.1, we have
$f(c^{\sqsubset^{F, \eta}}(m)) = f(\eta_X^{-1}(c^{\sqsubset}_{FX}(m))) \leq \eta_Y^{-1}((Ff)(c^{\sqsubset}_{FX}(\eta_X (m)))) \leq \eta_Y^{-1}(c^{\sqsubset}_{FX}((Ff)(\eta_X (m)))) \leq \eta_Y^{-1}(c^{\sqsubset}_{FX}(\eta_Y (f(m)))) \leq c^{\sqsubset^{F, \eta}}(f(m))$. Since $\eta_X (c^{\sqsubset^{F, \eta}}(m)) \leq c^{\sqsubset}_{FX}(\eta_X (m)), F \mbox{ is } (c^{\sqsubset^{F, \eta}}, c^{\sqsubset})$-continuous. If $c^{'}$ is another closure operator such that $(c^{'}; c^{\sqsubset})$-continuous, then $\eta_X (c^{'}(m)) \leq c^{\sqsubset}(\eta_X (m)) \Leftrightarrow c^{'}(m) \leq \eta_X^{-1}(c^{\sqsubset}_X (\eta_X (m))) = c^{\sqsubset^{F, \eta}}(m)$. If $\sqsubset$ is interpolative,then $c^{\sqsubset^{F, \eta}}(c^{\sqsubset^{F, \eta}}(m))\\  = c^{\sqsubset^{F, \eta}}(\eta_X^{-1}(c_{FX}^{\sqsubset}(\eta_X (m)))) = \eta_X^{-1}(c^{\sqsubset}_{FX}(\eta_X (\eta_X^{-1}(c^{\sqsubset}_{FX}(\eta_X (m)))))) \leq \eta_X^{-1}(c^{\sqsubset}_{FX}(\eta_X (m))) = \\ c^{\sqsubset^{F, \eta}}(m)$.
\end{proof}
\begin{Proposition}
	Let $\mathcal{A}$ be a reflective subcategory of $\mathcal{C}$ and $\sqsubset \in \bigwedge - TORD$. Then $c^{\mathcal{A}}(m) = \eta_X^{-1}(c^{\sqsubset}_{FX}(\eta_X (m)))$ is the largest closure operator on $\mathcal{C}$ for which every reflection morphism $\eta_X :
	X \longrightarrow FX$ is $(c^{\mathcal{A}}, c^{\sqsubset})$-continuous. Moreover, if $\sqsubset$ is interpolative, then $c^{\mathcal{A}}$ is idempotent.
\end{Proposition}
The above closure operator was first studied on the category of topological spaces and continuous maps by L. Stramaccia (\cite{stramaccia1988classes}), on topological categories by D. Dikranjan (\cite{dikranjan1992semiregular}) and later on an arbitrary category by D. Dikranjan and W. Tholen (\cite{MR2122803}). It is a special case of the pullback closure studied by D. Holgate in \cite{holgate1996pullback, holgate1995pullback} .
 \begin{Proposition}
 	Let $f^{-1}$ have a left adjoint for any $f \in \mathcal{C}$ and $(F, \eta)$ be a pointed endofunctor of $\mathcal{C}$ and $\sqsubset \in \bigvee - TORD$. Then $i^{\sqsubset^{F, \eta}}_X (m) = \eta_X^{-1}(i^{\sqsubset}_{FX}(\eta_X) (m))$ is the least interior operator on $\mathcal{C}$ for which every $\eta_X : X \longrightarrow FX$ is $(i^{\sqsubset}, i^{\sqsubset^{F, \eta}})$-continuous. $i^{F, \eta}$ is idempotent provided $\sqsubset$ is interpolative and each $\eta_X \in \varepsilon$.
 \end{Proposition}
 \begin{proof}
(I1) and (I2) is clearly satisfied. For (I3), let $f : X \longrightarrow Y$ be a $\mathcal{C}$-morphism and
 $m \in subY$. Using $i^{\sqsubset}$-continuity, Definition 2.1 and Lemma 2.1, we have, $f^{-1}(i^{\sqsubset^{F,\eta}}_Y (m)) = f^{-1}(\eta_X^{-1}(i^{\sqsubset}_{FX}(\eta_X)(m))) = \eta_X^{-1}(Ff)^{-1}(i^{\sqsubset}_{FY} (\eta_Y)(m)) \leq \eta_X^{-1}(i^{\sqsubset}_{FY}((Ff)^{-1}(\eta_Y)(m))) \\ \leq \eta_X i^{\sqsubset}_{FX}(\eta_X) (f^{-1}(m)) = i^{\sqsubset^{F, \eta}}_Y (f^{-1}(m))$. Since $\eta_X^{-1}(i^{\sqsubset}_{FX} (n)) \leq \eta_X^{-1}(i^{\sqsubset}_{FX}(\eta_X) \eta_X (n)) = \\ i^{\sqsubset^{F,\eta}}_X (\eta_X^{-1}(n)), \eta_X$ is $(i^{\sqsubset^{F,\eta}}, i^{\sqsubset})$-continuous for any $X \in \mathcal{C}$ and $n \in subX$. If $i^{'}$ is another interior operator on $\mathcal{C}$ such that $\eta_X$ is $(i^{'},i^{\sqsubset})$-continuous, then $i^{\sqsubset^{F,\eta}}_X (m) = \eta_X^{-1}(i^{\sqsubset}_{FX}((\eta_X) (m))) \leq i^{'}_X (\eta_X^{-1}(\eta_X) (n)) \leq i^{'}_X (n)$. If $\sqsubset \in \bigvee - INTORD(\mathcal{C},\mathcal{M})$ and $\eta_X \in \varepsilon$ for all $X \in \mathcal{C}$, then $i^{\sqsubset^{F,\eta}}(i^{\sqsubset^{F,\eta}}(m)) = i^{\sqsubset^{F,\eta}}(\eta_X^{-1}(i^{\sqsubset}_{FX}(\eta_X) (m))) = \eta_X^{-1}(i^{\sqsubset}_{FX}(\eta_X) (\eta_X^{-1}(i^{\sqsubset}_{FX}(\eta_X) (m)))) = \\ \eta_X^{-1}(i^{\sqsubset}_{FX}(i^{\sqsubset}_{FX}(\eta_X) (m))) = \eta_X^{-1}(i^{\sqsubset}_{FX}(\eta_X) (m)).$
  \end{proof}
 While the topogenous order induced by a pointed endofunctor (Theorem 5.2) was obtained with the help of $\sqsubset$-initial morphism, the notion of  $\sqsubset$-weakly final morphism is used in the next theorem to obtain the topogenous order induced by a co-pointed endofunctor.
 \begin{Theorem}
 Let $(G, \varepsilon)$ be a copointed endofunctor of $\mathcal{ C }$ and $\sqsubset$ a topogenous order on $\mathcal{ C }$,
 then for all $m \in subX$ and $n \ge m, m \sqsubset^{G, \varepsilon} n \Leftrightarrow \varepsilon^{-1}_{X}(n) \sqsubset_{GX}\varepsilon^{-1}_{X}(n)$ is a topogenous order on $\mathcal{C}$. It is the largest one for which every $\varepsilon_X : GX \longrightarrow X$ is $(\sqsubset, \sqsubset^{G,\varepsilon})$-continuous.
 \end{Theorem}
 \begin{proof}
 $(T1)$ and $(T2)$ are easily seen to be satisfied. For $(T3)$, let $f : X \longrightarrow Y$ be a Cmorphism. Then for all $m \in subX$ and $n \in subY$ such that $f(m) \leq n, f(m) \sqsubset^{G, \varepsilon} n \Leftrightarrow
 \varepsilon^{-1}_{Y}(f(m)) \sqsubset_{GY}\varepsilon^{-1}_{Y} 1(n) \Rightarrow (Gf)(\varepsilon^{-1}_{X}(m)) \sqsubset{GY} \varepsilon^{-1}_{X}(n) \Rightarrow \varepsilon^{-1}_{X}(m) \sqsubset_{GX} \varepsilon^{-1}_{X}(f^{-1}(n)) \Rightarrow \sqsubset^{GX}_{X}
 f^{-1}(n)$. Now, $\varepsilon_X : GX \longrightarrow X$ is trially $(\sqsubset,\sqsubset^{G, \varepsilon})$-continuous and if $\sqsubset^{'}$ is another topogenous order on $\mathcal{ C }$ such that $\varepsilon$ is $(\sqsubset, \sqsubset^{'})$-continuous, then $m\sqsubset^{'}_{X} n \Rightarrow \varepsilon_{X}(\varepsilon^{-1}_{X}(m)) \sqsubset^{'} n \Rightarrow \varepsilon^{-1}_{X} (m) \sqsubset_X \varepsilon^{-1}_{X} \Leftrightarrow m \sqsubseteq^{G,\varepsilon}_{X} n$ for all $n \ge m$.
  \end{proof}
Viewing a coreflector as a copointed endofunctor, we get the next proposition.
\begin{Proposition}
	Let $\mathcal{A}$ be a coreflective subcategory of $\mathcal{C}$ and $\sqsubset$ topogenous order on
	$\mathcal{A}$, then for all $m \in subX$ and $n \ge m, m \sqsubset^{\mathcal{A}} n \Leftrightarrow \varepsilon^{-1}_{X}(n) \sqsubset_{GX} \varepsilon^{-1}_{X}(n)$ is a topogenous
	order on $\mathcal{C}$. It is the largest one for which every coreflection morphism $\varepsilon_{X} : GX \longrightarrow X$ is
	$(\sqsubset, \sqsubset^{\mathcal{A}})$-continuous.
\end{Proposition}
\begin{Proposition}
	Let $(G, \varepsilon)$ be a copointed endofunctor of $\mathcal{A}$ and $\sqsubset \in \bigwedge-TORD$, then for all $m \in subX, c^{\sqsubset^{G,\varepsilon}}(m) = m \vee \varepsilon_{X}(c^{\sqsubset}(\varepsilon^{-1}_{X}(m)))$ is a closure oparator on $\mathcal{C}$. It is the least closure operator for which every $\varepsilon_{X} : GX \longrightarrow X$ is $(c, c^{G,\varepsilon})$-continuous.
\end{Proposition}
\begin{proof}
$(C1)$ and $(C2)$ are easily seen to be satisfied. To check $(C3)$, let $f : X \longrightarrow Y$ be a $\mathcal{C}$-morphism, and $m \in subX, f(c^{\sqsubset^{G,\varepsilon}}(m)) = f(m \vee \varepsilon_X(c^{\sqsubset}(\varepsilon^{-1}_{X}1(m)))) = f(m) \vee
f(\varepsilon_{X}(c^{\sqsubset}(\varepsilon^{-1}_{X}(m)))) = f(m) \vee \varepsilon_{Y}(Gf)(c^{\sqsubset}(\varepsilon^{-1}_{X}(m)))) \leq f(m) \vee \varepsilon_{Y} (c^{\sqsubset}(Gf)(\varepsilon^{-1}_{X}(m)))) \leq f(m)\vee
\varepsilon_{Y} (c^{\sqsubset}(\varepsilon^{-1}_{Y}f((m)))) = c^{\sqsubset^{G,\varepsilon}}f(m))$. For any $X \in \mathcal{C}, \varepsilon_{X}$ is $(c^{\sqsubset},c^{\sqsubset^{G,\varepsilon}})$-continuous \\ since $\varepsilon_{X}(c^{\sqsubset}(\varepsilon^{-1}_{X}(m)) \leq c^{\sqsubset^{G,\varepsilon}}(m) \Leftrightarrow c^{\sqsubset}(\varepsilon^{-1}_X(m)) \leq \eta^{-1}_{X}(c^{\sqsubset^{G,\varepsilon}}(m))$. If $c^{'}$ is another topogenous order on $\mathcal{C}$ such that $\varepsilon_X$ is $(c;, c^{'})$-continuous, then $c^{\sqsubset}(\varepsilon^{-1}_{X}(m))\leq \varepsilon^{-1}_{X}(c^{'}(m)) \Leftrightarrow \varepsilon_{X}(c^{\sqsubset}(\varepsilon^{-1}_X(m)) \leq
c^{'}(m) \Rightarrow m \vee \varepsilon_X(c^{\sqsubset}_{X}(m)) \leq c^{'}(m) \Rightarrow c^{\sqsubset^{G,\varepsilon}}(m)) \leq c^{'}(m)$. $\square $.
\end{proof}
\begin{Proposition}
Let $\mathcal{A}$ be a coreflective subcategory of $\mathcal{C}$ and $\sqsubset \in \bigwedge-TORD$, then for all $m \in subX, c^{\mathcal{A}}(m) = m \vee \varepsilon_{X}(c^{\sqsubset}(\varepsilon^{-1}_{X}(m)))$ is a closure operator on $\mathcal{C}$. It is the least closure
	operator for which every coreflection morphism $\varepsilon_{X} : GX \longrightarrow X$ is $(c, c^{\mathcal{A}})$-continuous.
\end{Proposition}
\begin{Proposition}
Assume that for any morphism $f \in \mathcal{C}$, the inverse image $f^{-1}$ commutes with the joins of subobjects. Let $(G, \varepsilon)$ be a copointed endofunctor of $\mathcal{C}$ and $\sqsubset \in  \bigvee-TORD$,
then for all $m \in subX, i^{\sqsubset^{G,\varepsilon}}(m) = m \wedge (\varepsilon_X)_{*}(i^{\sqsubset}(\varepsilon^{-1}_X(m)))$ is a interior operator on $\mathcal{C}$. It is the largest interior operator for which every $\varepsilon_X : GX \longrightarrow X$ is $(i, i^{G,\varepsilon})$-continuous.
\end{Proposition}
\begin{proof}
$(I1)$ and $(I2)$ are clearly seen to be satisfied. \\
For $(I3)$, let $f : X \longrightarrow Y$ be a $\mathcal{C}$-morphism, then $f^{-1}(i^{\sqsubset^{G,\varepsilon}}_Y(m)))= f^{-1}(m\wedge(\varepsilon_{Y})_{*}(i^{\sqsubset}(\varepsilon^{-1}_{Y}(m))))\\ = f^{-1}(m) \wedge f^{-1}((\varepsilon_{Y})_{*}(i^{\sqsubset}_{GY} (\varepsilon^{-1}_{Y}(m)))) \leq f^{-1}(m) \wedge (\varepsilon_{X})_{*}(Gf)^{-1}(i^{\sqsubset}_{GY} (\varepsilon^{-1}_{Y} (m)))) \leq f^{-1}(m) \wedge (\varepsilon_X)_{*}(i^{\sqsubset}_{GX}(Gf)^{-1}((\varepsilon^{-1}_{Y}(m)))) \leq f^{-1}(m) \wedge (\varepsilon_X)_{*}(i^{\sqsubset}_{GX}(\varepsilon^{-1}{X}(f^{-1}(m))) = i^{\sqsubset^{G, \varepsilon}}_{X}(f^{-1}(m)))$. \\
Since $\varepsilon^{-1}_{X}(i^{\sqsubset^{G,\varepsilon}}(m)) = \varepsilon^{-1}_{X}m)\wedge \varepsilon^{-1}_{X}((\varepsilon_X)_{*}(i^{\sqsubset}_{GX}(\varepsilon^{-1}X(m))) \leq i^{\sqsubset}_{GX}(\varepsilon^{-1}_{X}(m)) \varepsilon_X$ is $(i^{\sqsubset}, i^{\sqsubset^{G,\varepsilon}})$-continuous. If $i^{'}$ is another interior operator on $\mathcal{C}$ such that  $\varepsilon_{X}$ is $(i^{\sqsubset}, i^{'})$-continuous, then $\varepsilon^{-1}_{X}(i^{'}(m)) \leq i^{\sqsubset}(\varepsilon^{-1}_{GX}(n)) \Leftrightarrow i^{'}(m) \leq (\varepsilon_X)_{*}(i^{\sqsubset}_{GX}(\varepsilon^{-1}_{X}(n)) \Rightarrow i^{'}(n) \leq n\wedge (\varepsilon_X)_{*}(i^{\sqsubset}_{GX}(\varepsilon^{-1}_{X}(n))$.
\end{proof}
\section{Examples}

\begin{itemize}
\item [(1)] In the category $\bf Top $ of topological spaces and continuous maps with its (surjections, emdeddings)-factorization structure, consider the following two topogenous orders:
$A \sqsubset_X B \Leftrightarrow \bar{A} \subseteq B$ and $A \sqsubset^{'}_{X} B \Leftrightarrow A \subseteq B^{o}$ for all $X \in \bf Top $ and $A, B \subseteq X.$
Propositions $3.7$ and $3.8$ provide equivalent ways of caracterizing open and closed continuous maps as well as well as a continuous map whose domain carries the coarsest topology for which the map is continuous. Equivalent ways of describing hereditary quotient maps are provided. It is clear that $f_{*}(A) = Y\setminus f(X\setminus A)$ for any $A \subseteq X$.
\begin{Proposition}
A continuous map $f : X \longrightarrow Y$ is open if and only if $f$ is $\sqsubset$-co-strict if and only if $f$ is $\sqsubset^{'}$-strict.
\end{Proposition}
\begin{proof}
$ \overline{f^{-1}(A)} \subseteq B \Leftrightarrow f^{-1}\overline{(A)} \subseteq B \Leftrightarrow X \backslash B \subseteq X \backslash f^{-1}\overline{(A)} \Leftrightarrow  X \backslash B \subseteq f^{-1}(Y \backslash \bar{A}) \Leftrightarrow  f(X \backslash B) \subseteq Y \backslash \overline{A} \Leftrightarrow  \overline{A} \subseteq Y \backslash f(X \backslash B).$ $ B \subseteq (f^{-1}(A)]^{o} \Leftrightarrow  \overline{X\backslash f^{-1}(A)} \subseteq X \backslash B \Leftrightarrow  \overline{f^{-1}(Y\backslash A)} \subseteq X \backslash B \Leftrightarrow  \overline{Y\backslash A} \subseteq
Y \backslash f(B) \Leftrightarrow  f(B) \subseteq Y \backslash \overline{Y\backslash A} \Leftrightarrow  f(B) \subseteq A^{o}.$
$ B \subseteq (f^{-1}(A))^{o} \Leftrightarrow  f(B) \subseteq A^{o} \Leftrightarrow  B \subseteq f^{-1}(A^{o}).$
\end{proof}
A similar reasoning proves the following
\begin{Proposition}
A continuous map $f : X \longrightarrow Y $ is closed if and only if $f$ is $\sqsubset$-strict if and only if $f$ is $\sqsubset$-co-strict.
\end{Proposition}
For a continuous map $f : X \longrightarrow Y , X $ carries the initial topology induced by $f$ if,
$A \subseteq X$ is open if and only if there is an open $B \subseteq Y$ such that $A = f^{-1}(B)$ or equivalentely $\overline{A} = f^{-1}\overline{(f(A))}$ for each $A \subseteq X$ (see e.g \cite{ MR1039321, n1998general}). Such morphism corresponds to the $\sqsubset$-initial.
In fact, if $f$ is $\sqsubset$-initial, then $f^{-1}\overline{(f(A))} \subseteq B \Leftrightarrow X \backslash B \subseteq X \backslash f^{-1}(\overline{f(A)})  \Leftrightarrow f(X \backslash B) \subseteq Y \backslash \overline{f(A)}   \Leftrightarrow \overline{f(A)} \subseteq Y \backslash f(X \backslash B) \Leftrightarrow f(A) \sqsubset_Y f_{*}(B) \Leftrightarrow A \sqsubset B \Leftrightarrow \overline{A} \subseteq B $. Conversely if $X$ carries the coarsest topology for which $f$ is continuous,$A \sqsubseteq_{X} B \Leftrightarrow \overline{A} \subseteq B \Leftrightarrow f^{-1}(\overline{f(A)}) \subseteq B \Leftrightarrow X \backslash B \subseteq X \backslash f^{-1}(Y \backslash \overline{f(A))} \Leftrightarrow f(X \backslash B) \subseteq Y \backslash \overline{f(A)} \Leftrightarrow \overline{f(A)} \subseteq Y \backslash f(X \backslash B) \Leftrightarrow f(A) \sqsubset_Y f_{*}(B) $ for any $A, B \subseteq X$.  

An analogous reasoning shows that $X$ carries the coarsest topology for which f is continuous if and if $f(A) \subseteq [Y \backslash f(X \backslash B)]^{o} \Leftrightarrow A \subseteq B^{o}$. According to (\cite{MR1039321}, Exercise $2$F), $f$ is hereditary quotient if it is surjective with the property that every restriction $f^{'} : f^{-1}(A) \longrightarrow A$ is quotient for any $A \subseteq Y$ or equivalently $f$ surjective with the property that $f\overline{(f^{-1}(B))} \subseteq Y$ is closed for every $ B \subseteq Y$ . We clearly see that $f$ is $\sqsubset$-final if and only if it is $\sqsubset^{'}$-final if and only if it is hereditary quotient.

\item [(2)] In the category $\bf Grp $ of groups and group homomorphisms with the (epi, mono)-factorization structure, let  $\sqsubset$ be the topogenous structure defined by $A \sqsubset_G B \Leftrightarrow A \leq N \leq B$ with $ N \lhd G \mbox{and} A, B \leq G.$ A group homomorphism is  $\sqsubset$-strict if and only if it preserves normal subgroups. A group homomorphism $f : G \longrightarrow H $ is $\sqsubset$-final if and only if it is
surjective. Assume $f$ is $\sqsubset$-final.  Since  $H \lhd H$,  $f^{-1}(H) \lhd  G$.  Thus $f^{-1}(H)   \sqsubset_{G} G \leq f^{-1}(f(G)) \Rightarrow H \sqsubset_{H}  f(G)  \Rightarrow H \leq  f(G)  \Rightarrow   H = f(G) $, that is $f$ is surjective. Conversely if $f$ is surjective then it preserves normal subgroups and by Proposition 3.5(4), $f$ is $\sqsubset$-final. Proposition 3.5(3) allows to say that every injective group homomorphism that preserves normal subgroups is $\sqsubset$-initial. 

\item [(4)] Let $ \bf Prox$  be the category of proximity spaces and proximal maps with (surjective, embedding)-factorization. For any $(X, \delta) \in$ $ \bf Prox$  and $A, B \subseteq X, A \sqsubset_{(X, \delta)} B \Leftrightarrow A\bar{\delta}(X \backslash B)$ is an interpolative topogenous order on $ \bf Prox$.

\item [(5)] Consider $ \bf TopGrp$, the category of topological groups and continuous group homomorphisms. The forgetful functor $F$ : $\mbox{\textbf{TopGrp}} \longrightarrow \mbox{\textbf{Grp}}$ is a mono-fibration since every subgroup of a topological group is a topological group with the subspace topology (see [4]). Thus by Proposition 4.2, every topogenous order on $ \bf Grp$ can be initially lifted to a topogenous
order on $ \bf TopGrp$.

\item [(6)] Let $ \bf Top$ be the category of topological spaces and continuous maps with its (surjections, emdeddings)-factorization structure. It is well known that  $ \bf Top$, the category of  $T_o$-topological spaces and continuous maps is a epi-reflective subcategory of  $ \bf Top$. Define 
$S_X = \left\lbrace  \sqsubset_{X_{o}}  \vert X_{o} \in Top_{o} \right\rbrace$ 
by $ A \sqsubset_{X_{o}} B \Leftrightarrow \bar{A}\subseteq B $ for any $X_{o} \subseteq$ $Top_{o}, A, B \subseteq X_{o}$. Let $(F, \eta)$ be the reflector into  $ \bf Top$.  For any $X \in$  $ \bf Top$, $\eta_{X} : X \longrightarrow X/\sim$ takes each $x \in X$
to its equivalence class $[x] = \left\lbrace y \in X \, \vert \bar{\left\lbrace x\right\rbrace } = \bar{\left\lbrace y\right\rbrace }\right\rbrace $. Thus $S_X = \left\lbrace \sqsubset^{F,\eta}_{X} \, \vert X \in Top \right\rbrace$  with $ A \sqsubset^{F, \eta}_{X} B \Leftrightarrow \eta^{-1}_{X}\bar{\eta_{X}(A)} \subseteq B A, B \subseteq X$. 

\item [(6)] The category  $ \bf sTop$ of sequential topological spaces (those spaces in which every sequentially closed set is closed) is M-coreflective in $ \bf Top$. Consider $\sqsubset$ on $ \bf sTop$ defined by $A \sqsubset B \Leftrightarrow A \sqsubseteq C \sqsubseteq B$ for any $X \in $\textbf{sTop} and some closed subset $C$ of $X$. Let $(G, \varepsilon)$ be the coreflector into \textbf{sTop}. For any $(X, \mathcal{T^{'}} ) \in \textbf{Top}, \varepsilon_{(X,\mathcal{T^{'}})} : (X, \mathcal{T^{'}}) \longrightarrow (X, \mathcal{T} )$, identity map on $X$, where $\mathcal{T^{'}} = \left\lbrace A \sqsubseteq X  \vert  X  \backslash A \, is \, sequentially \, closed\,  in \, (X, \mathcal{T} ) \right\rbrace$  is an $ \bf sTop$-coreflection for any $(X, \mathcal{T} )$. It is clear that  $ A \sqsubset^{G,\varepsilon}_{X} B\Leftrightarrow A \subseteq C \subseteq B$ for some sequentially closed subset $C$ of $X$.
\end{itemize}

\bibliographystyle{abbrv}
\bibliography{references}
\addcontentsline{toc}{chapter}{References}
\end{document}